\documentclass{article}
 \textwidth 18 cm \textheight 10 in
\newtheorem{theorem}{Theorem}
\title{Counterexamples of the Geometrization Conjecture }
\author{Sze Kui Ng\\
{\small Department of Mathematics, Hong Kong Baptist University,
Hong Kong} \\{\small Email: szekuing@hotmail.com} }
\begin{document}
\date{}
\maketitle
\begin{abstract}
In this paper we propose  counterexamples to the Geometrization
Conjecture and the Elliptization Conjecture.

{\bf Mathematics Subject Classification: } 57M50, 57M27, 57N10,
57N16.
\end{abstract}

\section{A counterexample of the Geometrization Conjecture}\label{sec00}

A version of the Thurston's Geometrization Conjecture  states that
if a closed (oriented and connected) 3-manifold is irreducible and
atoroidal, then it is geometric in the sense that it can either
have a hyperbolic geometry or have a spherical geometry
\cite{Th1}\cite{Th2}\cite{Cas}. In this paper we propose
counterexamples to this conjecture by using the Dehn surgery
method of constructing closed 3-manifolds
 \cite{Ro}\cite{Li}.

 Let $K_{RT}^1$ denote the right trefoil knot with framing 1. Let
$K_{E}^r$ denote the figure-eight knot with framing $r$ where
$r=\frac pq$ is a rational number ($p$ and $q$ are co-prime
integers) such that $r>4$. We then consider a Dehn surgery on the
framed link $L=K_{RT}^1\cup K_{E}^r$ where the linking $\cup$ is
of the simplest Hopf link type.

We have that the Dehn surgery on $K_{RT}^1$ gives the Poincar{\'e}
sphere $M_{RT}^1$ which is with spherical geometry and with a
finite nontrivial fundamental group
\cite{Th1}\cite{Th2}\cite{Ro}\cite{Ht}\cite{Br}. Then
the Dehn surgery on $K_{E}^r$ gives a hyperbolic manifold
$M_{E}^r$ \cite{Th1}\cite{Th2}\cite{Ht}\cite{Br}. We want to show
that the 3-manifold $M_L$ obtained from surgery on $L$ is
irreducible and atoroidal, and is not geometric. From this we then
have that $M_L$ is a counterexample of the Geometrization
Conjecture.

Let us first show that $M_L$ is irreducible and atoroidal. From
\cite{Ng} we have the following quantum invariant
$\overline{W}(K_{RT}^1)$ of $M_{RT}^1$:
\begin{equation}
\overline{W}(K_{RT}^1)=R^{2}R_1^{-1}R_2^{1}W(C_1)W(C_2)
 \label{1}
\end{equation}
where the indexes of the $R$-matrices $R_1$ and $R_2$ are $1$ and
$-1$ respectively (These $R$-matrices are the monodromies of the
Knizhnik-Zamolodchikov equation; the notation $W(K)$ denotes the
generalized Wilson loop of a knot $K$ and is a quantum
representation of $K$ \cite{Ng}). Thus the indexes of
$R_1$ and $R_2$ are nonzero and are different. In \cite{Ng} we
call this property as the maximal non-degenerate property which is
a property only from nontrivial knots. We have that $R_1$ and
$R_2$ act on $W(C_1)$ and $W(C_2)$ respectively while $R$ is a
$R$-matrix for the linking of the framed knot $K_{RT}^1$ and acts
on $W(C_1)$ and $W(C_2)$. Similarly we have the following quantum
invariant of $M_{E}^r$:
\begin{equation}
\overline{W}(K_{E}^r)=R^{2p}R_1^{-3}R_2^{-a3}W(C_1)W(C_2)
 \label{2}
\end{equation}
where we choose a rational number $r=\frac pq$ such that the
integer $a\neq 1$ is nonzero. This is then
the maximal non-degenerate property.

Now let us consider the manifold $M_L$. Since $K_{RT}^1$ and
$K_{E}^r$ both have the maximal non-degenerate property we have
that there is no degenerate degree of freedom for the quantum
representation of $M_L$ by using the link $L$. From this we have
that $L$ is a minimal link for the Dehn surgeries obtaining $M_L$
\cite{Ng} (We shall later give more explanations on the definition
of minimal link and the related theorems on the classification of
3-manifolds by quantum invariant of 3-manifolds). It follows that
the quantum invariant of $M_L$ is given by the quantum
representation of $L$ and is of the following form:
\begin{equation}
\overline{W}(L)= P_L\overline{W}(K_{RT}^1)\overline{W}(K_{E}^r)
 \label{3}
\end{equation}
where $P_L$ denotes the linking part of the representation of $L$.

In this  quantum invariant (\ref{3}) of $M_L$ we have that
$\overline{W}(K_{RT}^1)$ and $\overline{W}(K_{E}^r)$ representing
$K_{RT}^1$ and $K_{E}^r$ respectively are independent of each
other and that the framed knots $K_{RT}^1$ and $K_{E}^r$ are
independent of each other in the sense that the framed knots
$K_{RT}^1$ and $K_{E}^r$ do not wind each other in the form as
described by the second Kirby move \cite{Ro}\cite{Ki}.

 We have that the
quantum invariant (\ref{3}) of $M_L$ uniquely represents $M_L$
because $L$ is minimal (We shall explain this point in the next
section). This means that there are no nontrivial symmetry
transforming it to another representation of $M_L$ with two framed
knots such that their quantum representations are different from
the two quantum representations $\overline{W}(K_{RT}^1)$ and
$\overline{W}(K_{E}^r)$ in (\ref{3}).

 Let us then first show that $M_L$ is irreducible. Since  the
quantum invariant (\ref{3}) of $M_L$ uniquely represents $M_L$ and
thus represents topological properties of $M_L$ we have that the
linking part $P_L$ of (\ref{3}) is a topological property of $M_L$
and thus cannot be eliminated. From this linking of
$\overline{W}(K_{RT}^1)$ and $\overline{W}(K_{E}^r)$ in (\ref{3})
we have that the invariant (\ref{3}) of $M_L$ cannot be written as
a free product form
$\overline{W}(K_1^{r_1})\overline{W}(K_2^{r_2})$ of two unlinked
framed knots $K_1^{r_1}$ and $K_2^{r_2}$ where each
$\overline{W}(K_i^{r_i}), i=1,2$ gives a closed 3-manifold. From
this we have that $M_L$ cannot be written as a connected sum of
two closed 3-manifolds. This shows that $M_L$ is irreducible.

Then we want to show that $M_L$ is atoroidal. Since the toroidal
property of a 3-manifold $M$ is about the existence of an infinite
cyclic subgroup $Z\oplus Z$ in $\pi_1(M)$ and is a property
derived from closed curves in $M$ only we have that this toroidal
property is derived from framed knots only since framed knots are
closed curves for constructing 3-manifolds. Now since $L$ is
minimal we have that the representation (\ref{3}) uniquely
represents $M_L$ and thus it gives all the topological properties
of $M_L$. From this we have that if $M_L$ has the toroidal
property then this property can only be derived from the two
framed knot components $K_{RT}^1$ and $K_{E}^r$. Now we have that
the 3-manifolds $M_{RT}^1$ and $M_{E}^r$ are both atoroidal and that 
the fundamental group of $M_{RT}^1$ is finite
\cite{Th1}\cite{Th2}\cite{Ht}\cite{Br}. Thus  the two
framed knot components $K_{RT}^1$ and $K_{E}^r$ do not give the
toroidal property of $M_L$. This shows that $M_L$ does not have the
toroidal property. Thus $M_L$ is atoroidal.

Let us explicitly compute the fundamental group
$\pi_1(M_L)$ of $M_L$ to give another proof for that $M_L$ is atoroidal. 
We have that $L=K_{RT}^1\cup K_{E}^r$ is of the Hopf link type. Thus by a computation similar to the computation of the link group of the Hopf link which is a direct product of the two knot groups of the two unknots forming the Hopf link we have that  the fundamental group $\pi_1(M_L)$ of $M_L$ is a direct product of the fundamental groups $\pi_1(M_{RT}^1)$ 
and $\pi_1(M_{E}^r)$:
\begin{equation}
\pi_1(M_L)=\pi_1(M_{RT}^1)\ast \pi_1(M_{E}^r)
 \label{3c}
\end{equation}
where $\pi_1(M_{RT}^1)\ast \pi_1(M_{E}^r)$ denotes the direct product of the fundamental groups $\pi_1(M_{RT}^1)$ and $\pi_1(M_{E}^r)$. 
Now since the 3-manifolds $M_{RT}^1$ and $M_{E}^r$ are both atoroidal and that the fundamental group $\pi_1(M_{RT}^1)$ is finite we have that
$\pi_1(M_L)$ does not contain a subgroup of the form $Z\oplus Z$.  This shows that $M_L$ does not have the
toroidal property. Thus $M_L$ is atoroidal.

Now since  the quantum invariant (\ref{3}) uniquely represents
$M_L$ we have that the two components $\overline{W}(K_{RT}^1)$ and
$\overline{W}(K_{E}^r)$ are topological properties of $M_L$. Then
since $\overline{W}(K_{RT}^1)$ (or $K_{RT}^1$) gives spherical
geometry property to $M_L$ and $\overline{W}(K_{E}^r)$ (or
$K_{E}^r$) gives hyperbolic geometry property to $M_L$ we have
that $M_L$ is not geometric. Indeed, since the two independent
components $\overline{W}(K_{RT}^1)$ and $\overline{W}(K_{E}^r)$ of
(\ref{3}) represent the manifolds $M_{RT}$ and $M_{E}$
respectively (and thus represent the fundamental groups
$\pi_1(M_{RT})$ and $\pi_1(M_{E})$ of $M_{RT}$ and $M_{E}$
respectively) we have that the fundamental group $\pi_1(M_L)$ of
$M_L$ contains the direct product $\pi_1(M_{RT})\ast\pi_1(M_{E})$
of the fundamental groups $\pi_1(M_{RT})$ and $\pi_1(M_{E})$. Now
let $\tilde{M_{L}}$ denote the universal covering space of $M_L$.
Then we have that $\pi_1(M_L)$ acts isometrically on
$\tilde{M_{L}}$. Now since $\pi_1(M_{RT})$ of the Poincar{\'e }
sphere $M_{RT}$ is not a subgroup of the isometry group of the
hyperbolic geometry $H^3$ and $\pi_1(M_{E})$ is not a subgroup of
the isometry group of the spherical geometry $S^3$ we have that
$\pi_1(M_{RT})\ast\pi_1(M_{E})$ is not a subgroup of the isometry
group of $H^3$ and is not a subgroup of the isometry group of
$S^3$. Thus $\pi_1(M_L)$ is not a subgroup of the isometry group
of $H^3$ and is not a subgroup of the isometry group of $S^3$. It
follows that $\tilde{M_{L}}$ is not the hyperbolic geometry $H^3$
and is not the spherical geometry $S^3$. This shows that $M_L$ is
not geometric, as was to be proved. Now since $M_L$ is irreducible
and atoroidal and is not geometric we have that $M_L$ is a
counterexample of the Geometrization Conjecture.

\section{Minimal link and classification of closed 3-manifolds}\label{sec01}

In this section we give more explanations on the definition of
minimal link and the related theorems on the classification of
closed 3-manifolds by quantum invariant used in the above
counterexample.

We have the following theorem of one-to-one representation of
3-manifolds obtained from framed knots $K^{\frac{p}{q}}$
\cite{Ng}:
\begin{theorem}
Let $M$ be a closed (oriented and connected) 3-manifold which is
constructed by a Dehn surgery on a framed knot $K^{\frac{p}{q}}$
 where $K$ is a nontrivial knot and $M$ is not a lens space. Then
 we have the following one-to-one representation of $M$:
\begin{equation}
\overline{W}(K^{\frac{p}{q}}):=R^{2p}R_1^{-m}R_2^{-am}W(C_1)W(C_2)
\label{qq1}
\end{equation}
where $m\neq 0$ ($m$ is also denoted by $m_1$ in \cite{Ng}) is the index of a nontrivial knot (which may or may not be the knot $K$ such that $M$ is also obtained from this
knot by Dehn surgery) and $am\neq 0$ is an integer related to
$m, p$ and $q$ such that $am\neq m$ (Thus (\ref{qq1}) is with
the maximal non-degenerate property).
\end{theorem}

We remark that if $M$ is a lens space we can also define a similar
quantum invariant $\overline{W}(K^{\frac{p}{q}})$ for
$M$ which however is not of the above maximal non-degenerate form
\cite{Ng}.

 Let us then consider a
3-manifold $M$ which is obtained from a framed link $L$ with the
minimal number $n$ of component knots where $n\geq 2$ (where the
minimal number $n$ means that if $M$ can also be obtained from
another framed link then the number of component knots of this
framed link must be $\geq n$) . In this case we call $L$ a minimal
link of $M$. From the generalized second Kirby moves (which generalizes second Kirby move from integer to rational number \cite{Ng} and for simplicity we shall call them again as the second kirby moves) 
 we may suppose that $L$ is
in the form that the components $K_i^{\frac{p_i}{q_i}},i=1,...,n$
of $L$ do not wind each other in the form described by the second
Kirby move. In this case we say that this minimal $L$ is in the
form of maximal non-degenerate state where the degenerate property
is from the winding of one component knot with the other component
knot by the second Kirby moves. Thus this $L$ has both the minimal
and maximal property as described. Then we want to find a
one-to-one representation (or invariant) of $M$ from this $L$. Let
us write $W(L)$, the generalized Wilson loop of $L$, in the
following form \cite{Ng}:
\begin{equation}
W(L)=P_L \prod_i W(K_i^{\frac{p_i}{q_i}}) \label{qq2a}
\end{equation}
where $P_L$ denotes a product of $R$-matrices acting on a subset
of $\{W(K_i),W(K_{ic}), i =1,...,n \}$ where
$W(K_i^{\frac{p_i}{q_i}})$ are independent (This is from the form
of $L$ that the component knots $K_i$ are independent in the sense
that they do not wind each other by the second Kirby moves). Then
we consider the following representation (or invariant) of $M$:
\begin{equation}
\overline{W}(L):=P_L \prod_i \overline{W}(K_i^{\frac{p_i}{q_i}})
\label{qq2}
\end{equation}
where we define $\overline{W}(K_i^{\frac{p_i}{q_i}})$ by
(\ref{qq1}) and they are independent. We then have  the following
theorem:
\begin{theorem}
Let $M$ be a closed (oriented and connected) 3-manifold which is
constructed by a Dehn surgery on a minimal link $L$ with the
minimal number $n$ of component knots (and with the maximal
property). Then we have that (\ref{qq2}) is
 a one-to-one representation (or invariant) of $M$.
\end{theorem}

{\bf Proof}. We want to show
that (\ref{qq2}) is a one-to-one representation (or invariant) of
$M$. Let $L^{\prime}$ be another framed link for $M$ which is also
with the minimal number $n$ (and with the maximal property). Then
we want to show $\overline{W}(L)=\overline{W}(L^{\prime})$.

For the case $n=1$ this is true by the above theorem for manifolds $M$ obtained from minimal framed knot $K^{\frac{p}{q}}$.

Let us consider $n\geq 2$. 
Since the components of $L$ do not wind each other as described by the second Kirby move we have that the components of $L$ are independent of each other. Thus there is no nontrivial homeomorphism changing these components $\overline{W}(K_i^{\frac{p_i}{q_i}})$  except those homeomorphisms involving the second Kirby moves for the winding of the components of $L$ with each other. 
Then under  the second Kirby moves we have that the components of $L$ wind each other and thus will reduce the independent degree of freedom to be less than $n$. Thus to restore the degree of freedom to $n$ these homeomorphisms must also contain the first Kirby moves of adding unknots with framing $\pm 1$. In this case these unknots can be deleted and thus $L$ 
 is not minimal and this is a contradiction. Thus there is no nontrivial homeomorphism changing the components $\overline{W}(K_i^{\frac{p_i}{q_i}})$ of $\overline{W}(L)$ except those homeomorphisms consist of only the second Kirby moves for the winding of the components of $L$ with each other. 

 Now suppose that $\overline{W}(L)\neq \overline{W}(L^{\prime})$. Then there exists nontrivial homeomorphism of changing $L$ to $L^{\prime}$ for changing the components $\overline{W}(K_i^{\frac{p_i}{q_i}})$ of $\overline{W}(L)$ to the components of 
 $\overline{W}(L^{\prime})$. This is impossible since there are no  nontrivial homeomorphsm for changing these components $\overline{W}(K_i^{\frac{p_i}{q_i}})$ except those homeomorphisms consist of only the second Kirby moves for the winding of the components of $L$ with each other. Thus $\overline{W}(L)=\overline{W}(L^{\prime})$.

Thus we have that (\ref{qq2}) is a one-to-one representation (or invariant) of $M$, as was to be proved. 
$\diamond$

As a converse to the above theorem let us suppose that  the
representation (\ref{qq2}) uniquely represents $M_L$ in the sense
that there are no nontrivial symmetry  
transforming the $n$
independent components of $\overline{W}(L)$ to other $n$
independent components of $\overline{W}(L^{\prime})$ where the
link $L^{\prime}$ also gives the manifold $M_L$. Then from the
above proof we see that the link $L$ is a minimal (and maximal) link for
obtaining $M_L$.

{\bf Remark}. Let $L$ be a minimal (and maximal) framed link. Then from the above proof we have that the components of $L$ are independent of each other in the sense that if we transform a component framed knot of $L$ to an equivalent framed knot by a homeomorphism then the other components of $L$ are not affected by this transformation. $\diamond$

Now let us consider the framed link $L=K_{RT}^1\cup K_{E}^r$ in
the above section. We have that the knot components $K_{RT}^1$ and
$K_{E}^r$ of $L$ do not wind each other in the form as described
by the second Kirby move. Thus we have that their corresponding
quantum invariants $\overline{W}(K_{RT}^1)$ and
$\overline{W}(K_{E}^r)$ are independent. Then
$\overline{W}(K_{RT}^1)$ and $\overline{W}(K_{E}^r)$ are in the maximal non-degenerate
form which is invariant under all homeomorphisms execept the second Kirby moves which are excluded (Indeed for $\overline{W}(K_{RT}^1)$ there is a homeomorphism transforming $K_{RT}^1$ to $K_{E}^{-1}$. Then the informations of these two frame knots are included in $\overline{W}(K_{RT}^1)$ and thus $\overline{W}(K_{RT}^1)$ is invariant under this homeomorphism. Then since $\overline{W}(K_{RT}^1)$ is in the maximal non-degenerate
form there are no degenerate degree of freedoms for other homeomorphisms execept the second Kirby moves which reduce the degree of freedom of $L$. Similarly for $\overline{W}(K_{E}^r)$). Thus $L$ is a minimal (and maximal) link of $M_L$ and the representation (\ref{3}) is the quantum invariant of $M_L$.

\section{A counterexample of the Elliptization Conjecture}\label{sec03}

The above counterexample of the Geometrization Conjecture is with
an infinite fundamental group. Let us in this section propose a
counterexample which is with a finite fundamental group to the
Geometrization Conjecture. This example is then also a
counterexample of the Thurston's Elliptization Conjecture which
states that if a closed (oriented and connected) 3-manifold is
irreducible and atoroidal and is with a finite fundamental group
then it is geometric in the sense that it can have a spherical
geometry \cite{Th1}\cite{Th2}\cite{Cas}.

Let us consider a Dehn surgery on the framed link $L=K_{RT}^1\cup
K_{RT}^1$ where the linking $\cup$ is of the simplest Hopf link
type. We want to show that the 3-manifold $M_L$ obtained from this
surgery is a counterexample of the Elliptization Conjecture.

As similar to the above example we have that this $L$ is minimal
and the 3-manifold $M_L$ is uniquely represented by the following
quantum invariant:
\begin{equation}
\overline{W}(L)= P_L\overline{W}(K_{RT}^1)\overline{W}(K_{RT}^1)
 \label{3a}
\end{equation}
where $P_L$ denotes the linking part of the representation of $L$.

Then as similar to the above example we have that this 3-manifold
$M_L$ is irreducible and atoroidal. Let us then show that $M_L$ is
with a finite fundamental group and is not geometric. Since  the
quantum invariant (\ref{3a}) uniquely represents $M_L$ we have
that the two components $\overline{W}(K_{RT}^1)$ are topological
properties of $M_L$. Then we have that the fundamental group
$\pi_1(M_L)$ of $M_L$ contains the direct product
$\pi_1(M_{RT})\ast\pi_1(M_{RT})$.

Further as similar to the above example
because $L$ is of the Hopf link type we have that $\pi_1(M_L)=\pi_1(M_{RT}^1)\ast \pi_1(M_{RT}^1)$. Now since the fundamental group $\pi_1(M_{RT})$ is finite 
we have that the fundamental group $\pi_1(M_L)$ is also finite.

  Now let
$\tilde{M_{L}}$ denote the universal covering space of $M_L$. Then
we have that $\pi_1(M_L)$ acts isometrically on $\tilde{M_{L}}$.
We want to show that $\tilde{M_{L}}$ is not the 3-sphere $S^3$.
Suppose this is not true. Then since $\pi_1(M_L)$ contains (and equals to) the
direct product $\pi_1(M_{RT})\ast\pi_1(M_{RT})$ we have that the
direct product $\pi_1(M_{RT})\ast\pi_1(M_{RT})$ is a subgroup of
the isometry group of $S^3$. Now since $S^3$ is a fully isotropic
manifold containing no boundary ($S^3$ is closed) there is no way
to distinguish two identical but independent subgroups
$\pi_1(M_{RT})$ of the isometry group of $S^3$. From this we have
that the direct product $\pi_1(M_{RT})\ast\pi_1(M_{RT})$ can only
act on $S^3\times S^3$ where each $\pi_1(M_{RT})$ acts on a
different $S^3$ and cannot act on the same $S^3$ such that
 $\pi_1(M_{RT})\ast\pi_1(M_{RT})$ acts on
 $S^3$ (Comparing to the hyperbolic case we have that the
direct product of two subgroups of the isometry group of the
hyperbolic geometry $H^3$ may act on $H^3$ since $H^3$ has
nonempty boundary which can be used to distinguish two identical
but independent subgroups of the isometry group of $H^3$). Thus
the direct product $\pi_1(M_{RT})\ast\pi_1(M_{RT})$ is not a
subgroup of the isometry group of $S^3$ (We can also prove this
statement by the fact that $\pi_1(M_{RT})$ is a nonabelian
subgroup of the rotation group $O(4)$ which is the isometry group
of $S^3$. Indeed since $\pi_1(M_{RT})$ is nonabelian it must act
on a space with dimension $\geq 3$. Thus
$\pi_1(M_{RT})\ast\pi_1(M_{RT})$ must act on a space with
dimension $\geq 6$. Now $O(4)$ can only act on a space with
dimension $4$ we have that $\pi_1(M_{RT})\ast\pi_1(M_{RT})$ is not
a subgroup of $O(4)$). This is a contradiction. This contradiction
shows that $\tilde{M_{L}}$ is not the 3-sphere $S^3$. Thus $M_L$
is not geometric.
 Now since $M_L$ is irreducible
and atoroidal and is with finite fundamental group and is not
geometric we have that $M_L$ is a counterexample of the
Elliptization Conjecture.

 \end{document}